%

\catcode`\@=11

\magnification=1200
\baselineskip=14pt

\pretolerance=500    \tolerance=1000 \brokenpenalty=5000

\catcode`\;=\active
\def;{\relax\ifhmode\ifdim\lastskip>\z@
\unskip\fi\kern.2em\fi\string;}

\overfullrule=0mm

\catcode`\!=\active
\def!{\relax\ifhmode\ifdim\lastskip>\z@
\unskip\fi\kern.2em\fi\string!}

\catcode`\?=\active
\def?{\relax\ifhmode\ifdim\lastskip>\z@
\unskip\fi\kern.2em\fi\string?}

\frenchspacing

\newif\ifpagetitre            \pagetitretrue
\newtoks\hautpagetitre        \hautpagetitre={ }
\newtoks\baspagetitre         \baspagetitre={1}

\newtoks\auteurcourant        \auteurcourant={ }
\newtoks\titrecourant
\titrecourant={ }
\newtoks\hautpagegauche       \newtoks\hautpagedroite
\hautpagegauche={\hfill\sevenrm\the\auteurcourant\hfill}
\hautpagedroite={\hfill\sevenrm\the\titrecourant\hfill}

\newtoks\baspagegauche       \baspagegauche={\hfill\rm\folio\hfill}

\newtoks\baspagedroite       \baspagedroite={\hfill\rm\folio\hfill}

\headline={
\ifpagetitre\the\hautpagetitre
\global\pagetitrefalse
\else\ifodd\pageno\the\hautpagedroite
\else\the\hautpagegauche\fi\fi}

\footline={\ifpagetitre\the\baspagetitre
\global\pagetitrefalse
\else\ifodd\pageno\the\baspagedroite
\else\the\baspagegauche\fi\fi}

\def\date{\ {\the\day}\
\ifcase\month\or Janvier\or F\'evrier\or Mars\or Avril
\or Mai \or Juin\or Juillet\or Ao\^ut\or Septembre
\or Octobre\or Novembre\or D\'ecembre\fi\
{\the\year}}

\def\up#1{\raise 1ex\hbox{\sevenrm#1}}

\def\cqfd{\unskip\kern 6pt\penalty 500
\raise -2pt\hbox{\vrule\vbox to 10pt{\hrule width 4pt
\vfill\hrule}\vrule}\par\medskip}

\def\section#1{\vskip 7mm plus 20mm minus 1.5mm\penalty-50
\vskip 0mm plus -20mm minus 1.5mm\penalty-50
{\bf\noindent#1}\nobreak\smallskip}

\def\subsection#1{\medskip{\bf#1}\nobreak\smallskip}

\def\displaylinesno #1{\dspl@y\halign{
\hbox to\displaywidth{$\@lign\hfil\displaystyle##\hfil$}&
\llap{$##$}\crcr#1\crcr}}

\def\ldisplaylinesno #1{\dspl@y\halign{
\hbox to\displaywidth{$\@lign\hfil\displaystyle##\hfil$}&
\kern-\displaywidth\rlap{$##$}
\tabskip\displaywidth\crcr#1\crcr}}

\def\hfl#1#2{\smash{\mathop{\hbox to 12 mm{\rightarrowfill}}
\limits^{\scriptstyle#1}_{\scriptstyle#2}}}

%
%
\def\bibn@me{R\'ef\'erences}
\def\bibliographym@rk{\centerline{{\sc\bibn@me}}
	\sectionmark\section{\ignorespaces}{\unskip\bibn@me}
	\bigbreak\bgroup
	\ifx\ninepoint\undefined\relax\else\ninepoint\fi}
%
%
%
\let\refsp@ce=\
\let\bibleftm@rk=[
\let\bibrightm@rk=]
%
%
%
\def\numero{n\raise.82ex\hbox{$\fam0\scriptscriptstyle
o$}~\ignorespaces}
%
%
\newcount\equationc@unt
\newcount\bibc@unt
\newif\ifref@changes\ref@changesfalse
\newif\ifpageref@changes\ref@changesfalse
\newif\ifbib@changes\bib@changesfalse
\newif\ifref@undefined\ref@undefinedfalse
\newif\ifpageref@undefined\ref@undefinedfalse
\newif\ifbib@undefined\bib@undefinedfalse
\newwrite\@auxout
%
%
%
%
%
%
%
%
\def\re@dreferences#1#2{{%
	\re@dreferenceslist{#1}#2,\undefined\@@}}
\def\re@dreferenceslist#1#2,#3\@@{\def\next{#2}%
	\expandafter\ifx\csname#1@@\meaning\next\endcsname\relax
	??\immediate\write16
	{Warning, #1-reference "\next" on page \the\pageno\space
	is undefined.}%
	\global\csname#1@undefinedtrue\endcsname
	\else\csname#1@@\meaning\next\endcsname\fi
	\ifx#3\undefined\relax
	\else,\refsp@ce\re@dreferenceslist{#1}#3\@@\fi}
%
%
%
\def\newlabel#1#2{{\def\next{#1}\newl@bel#2}}
\def\newl@bel#1#2{%
	\expandafter\xdef\csname ref@@\meaning\next\endcsname{#1}%
	\expandafter\xdef\csname pageref@@\meaning\next\endcsname{#2}}
\def\label#1{{%
	\toks0={#1}\message{ref(\lastref) \the\toks0,}%
	\ignorespaces\immediate\write\@auxout%
	{\noexpand\newlabel{\the\toks0}{{\lastref}{\the\pageno}}}%
	\def\next{#1}%
	\expandafter\ifx\csname ref@@\meaning\next\endcsname\lastref%
	\else\global\ref@changestrue\fi%
	\newlabel{#1}{{\lastref}{\the\pageno}}}}
\def\ref#1{\re@dreferences{ref}{#1}}
\def\pageref#1{\re@dreferences{pageref}{#1}}
%
%
\def\bibcite#1#2{{\def\next{#1}%
	\expandafter\xdef\csname bib@@\meaning\next\endcsname{#2}}}
\def\cite#1{\bibleftm@rk\re@dreferences{bib}{#1}\bibrightm@rk}
%
%
\def\beginthebibliography#1{\bibliographym@rk
	\setbox0\hbox{\bibleftm@rk#1\bibrightm@rk\enspace}
	\parindent=\wd0
	\global\bibc@unt=0
	\def\bibitem##1{\global\advance\bibc@unt by 1
		\edef\lastref{\number\bibc@unt}
		{\toks0={##1}
		\message{bib[\lastref] \the\toks0,}%
		\immediate\write\@auxout
		{\noexpand\bibcite{\the\toks0}{\lastref}}}
		\def\next{##1}%
		\expandafter\ifx
		\csname bib@@\meaning\next\endcsname\lastref
		\else\global\bib@changestrue\fi%
		\bibcite{##1}{\lastref}
		\medbreak
		\item{\hfill\bibleftm@rk\lastref\bibrightm@rk}%
		}
	}
\def\endthebibliography{\egroup\par}
%
%
    \outer\def\bye{\@closeaux
    	\par\vfill\supereject\end}
%
\def\@closeaux{\closeout\@auxout
	\ifref@changes\immediate\write16
	{Warning, changes in references.}\fi
	\ifpageref@changes\immediate\write16
	{Warning, changes in page references.}\fi
	\ifbib@changes\immediate\write16
	{Warning, changes in bibliography.}\fi
	\ifref@undefined\immediate\write16
	{Warning, references undefined.}\fi
	\ifpageref@undefined\immediate\write16
	{Warning, page references undefined.}\fi
	\ifbib@undefined\immediate\write16
	{Warning, citations undefined.}\fi}
%
%
\immediate\openin\@auxout=\jobname.aux
\ifeof\@auxout \immediate\write16
     {Creating file \jobname.aux}
\immediate\closein\@auxout
\immediate\openout\@auxout=\jobname.aux
\immediate\write\@auxout {\relax}%
\immediate\closeout\@auxout
\else\immediate\closein\@auxout\fi
%
%
\input\jobname.aux \par
\immediate\openout\@auxout=\jobname.aux
%
%

\def\Z{{\bf Z}}
 
\def\R{{\bf R}}

\def\X{{\bf X}}

\def\be{{\bf e}}

\def\bP{{\bf P}}
\def\bQ{{\bf Q}}
\def\bR{{\bf R}}

\def\bx{{\bf x}}
\def\bX{{\bf X}}
\def\by{{\bf y}}
\def\bY{{\bf Y}}
\def\bz{{\bf z}}
\def\bZ{{\bf Z}}

\def\rem{\noindent {\bf{Remark. }}}  \def\pro{\noindent {\bf{Proof.
}}}

   \def\La{{\Lambda}}    
  \def\th{{\theta}} 
\def\Th{{\Theta}}

\def\om{{\omega}}
\def\omc{{\hat{\omega}}}
\def\and{\quad\hbox{and}\quad}
\def \lrcorner
{\mathrel{\hbox{\vrule depth 0 pt height 0.4pt width 4pt
\vrule depth 0 pt height 5pt width 0.4pt
\kern 1pt }}}

\def\build#1_#2^#3{\mathrel{\mathop{\kern 0pt#1}\limits_{#2}^{#3}}}



\def\pro{\noindent {\bf Proof. }}

\def\smallsquare{\vbox{\hrule\hbox{\vrule height 1 ex\kern 1
ex\vrule}\hrule}}
\def\cqfd{\hfill \smallsquare\vskip 3mm}

\def\hw{{\hat w}}

\def\cC{{\cal C}}
\def\cP{{\cal P}}
\def\vol{{\rm vol}}

\def\bibn@me{R\'ef\'erences bibliographiques}
%
\def\bibliographym@rk{\bgroup}
%
%
\outer\def\bye{ 	\par\vfill\supereject\end}


\null


\vskip 2mm

\centerline{\bf On transfer inequalities in Diophantine approximation, II}

\vskip 8mm

\centerline{Yann B{\sevenrm UGEAUD} \footnote{}{\rm
}
\& Michel L{\sevenrm AURENT}
}

\vskip 11mm

\section{ 1. Introduction}

 Let $n$ be a positive  integer 
and let $\Th =(\th_1, \dots , \th_n)$ be a point in $\bR^n$.
We shall assume in all the forthcoming statements 
that the real numbers $1, \th_1, \dots , \th_n$ 
are linearly independent over 
the field $\bQ$ of rational numbers. 
Khintchine's transference principle relates 
the sharpness of the rational simultaneous 
approximation to $\th_1, \dots , \th_n$
with the measure of linear independence over 
$\bQ$ of $1,\th_1, \dots , \th_n$. Let us first quantify these notions
by introducing the exponents $\om_0(\Th)$ and $\om_{n-1}(\Th)$
(the meaning  of the indices $0$ and $n-1$ will be explained
afterwards).

\proclaim Definition 1.
We denote respectively by $\om_0(\Th)$ and $\om_{n-1}(\Th)$ 
the supremum, possibly infinite,  
of the real numbers $\om$ for which there exist infinitely many 
integer $(n+1)$-tuples $(x_0, \dots , x_n)$ 
satisfying respectively the inequation
$$
\max_{1\le i\le n} |x_0 \th_i - x_i |  \le 
\Big( \max_{0\le i\le n} |  x_i | \Big)^{-\om}
\quad {\sl or} \quad 
|x_0  + x_1 \th_1 + \cdots +  x_n \th_n |  
\le \Big( \max_{0\le i\le n} |  x_i  | \Big)^{-\om}.
$$ 

Now we can state Khintchine's transference principle \cite{Kh26} (see
\cite{Mah36} for an alternative proof, and the
monographs \cite{Cas,SchmLN,GrLe}) as follows:

\proclaim 
Theorem K. The inequalities 
$$
{\om_{n-1}(\Th)\over (n-1)\om_{n-1}(\Th) +n} 
\le \om_0(\Th) \le { \om_{n-1}(\Th) -n+1\over n}
\leqno{(1.1)}
$$
hold for any point $\Th  =(\th_1, \dots , \th_n)$ in $\bR^n$ with 
$1, \th_1, \dots , \th_n$ linearly independent over $\bQ$.

Moreover, Jarn\'\i k \cite{Jar35,Jar36}
established that both inequalities in (1.1) are optimal, and, consequently,
that Theorem K is best possible. 
It is the main purpose of the present paper to show that,
however, Theorem K can be refined if we introduce two further
quantities associated with $\Th$.

Following the general ``hat'' notations of \cite{BuLa05a}, 
let us introduce the uniform
analogues of the exponents $\om_0(\Th)$ and $\om_{n-1}(\Th)$.

\proclaim Definition 2. 
We denote respectively by $\omc_0(\Th)$ and $\omc_{n-1}(\Th)$ 
the supremum of the real numbers $\om$ 
such that for all sufficiently large real number $X$,  
there exists a non-zero
integer $(n+1)$-tuples $(x_0, \dots , x_n)$ with supremum norm
$$
\max_{0\le i \le n} | x_i | \le X,
$$
satisfying respectively the inequation
$$
\max_{1\le i\le n} |x_0 \th_i - x_i |  \le X^{-\om}
\quad {\sl or} \quad 
 |x_0  + x_1 \th_1 + \cdots +  x_n \th_n |  \le X^{-\om}.
$$ 

We establish the 
following refinement of Khintchine's theorem, which 
involves the uniform exponents associated with $\Th$.

\proclaim 
Theorem 1. Suppose $n\ge 2$. The inequalities 
$$
{ (\omc_{n-1}(\Th) -1)\om_{n-1}(\Th) \over 
((n-2)\omc_{n-1}(\Th) +1)\om_{n-1}(\Th) +(n-1)\omc_{n-1}(\Th)}
\le \om_0(\Th)
$$
and
$$
\om_0(\Th) \le {(1 - \omc_0(\Th))\om_{n-1}(\Th) -n +2 -\omc_0(\Th) \over n-1}
$$
hold for any point $\Th  =(\th_1, \dots , \th_n)$ in $\bR^n$ with 
$1, \th_1, \dots , \th_n$ linearly independent over $\bQ$.

The above inequalities are stronger than (1.1), since 
$$
\omc_{n-1}(\Th) \ge n \and \omc_0(\Th) \ge {1\over n},
$$
by the Dirichlet Box Principle. 
Theorem 1 was first established  when $n=2$ in \cite{Lau08a}
and its statement was announced in \cite{BuLa07} and in \cite{Lau08b}.  It 
follows from the description given in \cite{Lau08a} 
of the set of all possible quadruples
$$
\bigl(\om_1(\Th),\om_0(\Th),\omc_1(\Th),\omc_0(\Th)\bigr), 
$$
where $\Th$ ranges over $\bR^2$,  that Theorem 1 is optimal in dimension two.

Theorem K was extended by Dyson \cite{Dys} to transfer inequalities
between approximation to a system of linear forms and
approximation of the tranpose system. It would be interesting to
establish a suitable extension of Theorem 1.

The present paper is organized as follows. 
In Section 2, we define further  exponents $\om_d (\Th)$ for  
$d = 1, \ldots , n-2$,
measuring the accuracy with which $\Th$ can be 
approximated by rational linear subvarieties of dimension $d$.
We state in Theorems 2 and 3 
transference inequalities linking $\om_d (\Th)$
and $\om_{d+1} (\Th)$, the composition of which gives
Theorem K. This was already known \cite{Schm67,Lau08b}, but our
proof, based on the second theorem of Minkowski, is new.
Furthermore, our method allows us to refine inequalities 
between $\om_0 (\Th)$ and $\om_1 (\Th)$
(resp. between $\om_{n-1} (\Th)$ and $\om_{n-2} (\Th)$),
by taking also $\omc_0 (\Th)$ (resp. $\omc_{n-1} (\Th)$)
into account. Using this, we get Theorem 1, as is
explained in Section 7.
Section 3 is devoted to some preliminaries
of multilinear algebra.
In Section 4 and at the beginning of Section 6, we give
alternative definitions of the exponents $\om_d$.
Theorems 2 and 3 are established in Sections 5 and 6, respectively.

\section{2. Going-up and going-down transfers}
It is convenient to view $\bR^n$ as a subset of $\bP^n(\bR)$
via the usual embedding $(x_1,\dots ,x_n)\mapsto (1,x_1,\dots ,x_n)$. 
We shall  identify $\Th  =(\th_1, \dots , \th_n)$
with its image in $\bP^n(\bR)$.

 Following  \cite{Lau08b}, 
let us introduce for each integer $d$ with 
$0 \le d \le n-1$ an exponent $\om_d(\Th)$ which measures  
the approximation to the point $\Th \in \bP^n(\bR)$
by rational linear projective subvarieties of dimension $d$, 
in terms of their height.  Denote by d the projective distance 
on $\bP^n(\bR)$ (it will be defined in \S 4 below; 
notice  however that the normalization used 
there does not matter for our purpose).
For any real linear subvariety $L$ of $\bP^n(\bR)$, we denote by 
$$
{\rm d} ( \Th , L) = \min_{ P\in L} {\rm d} (\Th , P)
$$
the minimal distance between $\Th$ and the real  points $P$ of $L$. 
When $L$ is rational over $\bQ$, we indicate  
moreover by $H(L)$ its height, that is the Weil height 
of any system of Pl\" ucker coordinates 
of $L$. It is  convenient to normalize the height 
by using the Euclidean norm at the Archimedean place of $\bQ$.
 We refer to  \S 1 of \cite{Schm67} for more 
information on the notion of height of a linear subspace.

\proclaim
Definition 3. Let $d$ be an integer with $0\le d \le n-1$. 
We denote by  $\om_d(\Th)$   the supremum  
of the real numbers $\om$ for which there exist 
infinitely many rational linear subvarieties $L \subset \bP^n(\bR)$
such that 
$$
\dim (L) = d \and {\rm d}(\Th , L) \le H(L)^{-1-\om}.
$$

Definitions 1 and 3 are consistent, since ${\rm d}(\Th, L)$ compares respectively with 
$$
\max_{1\le i\le n} \Big| \th_i -{x_i\over x_0}\Big|  
\and {| y_0 + y_1 \th_1 + \cdots + y_n \th_n| \over {\displaystyle\max_{0 \le i \le n} | y_i|}}
$$
when $L$ is either the rational point (case $d=0$) 
with homogeneous  coordinates $(1,x_1/x_0, \dots , x_n/x_0)$,  or the hyperplane 
(when $d=n-1$) with homogeneous equation $y_0 X_0 + \cdots + y_n X_n=0$.

Theorem 1 is a consequence of the following two
statements.

\proclaim
Theorem 2 {\bf (Going-up transfer)}. 
Let $\Th  =(\th_1, \dots , \th_n)$ be in $\bR^n$ with 
$1, \th_1, \dots , \th_n$ linearly independent over $\bQ$.
For any integer $d$ with $0\le d \le n-2$, we have the lower bound
$$
\om_{d+1}(\Th) \ge { (n-d) \om_d(\Th) +1 \over n-d-1}.
\leqno{(2.1)}
$$
Furthermore,
$$
\om_1(\Th) \ge {\om_0(\Th) + \omc_0(\Th) \over 1 - \omc_0(\Th)}. \leqno (2.2)
$$

\proclaim
Theorem 3  {\bf (Going-down transfer)}. 
Let $\Th  =(\th_1, \dots , \th_n)$ be in $\bR^n$ with 
$1, \th_1, \dots , \th_n$ linearly independent over $\bQ$.
For any integer $d$ with $1\le d \le n-1$, we have the lower bound
$$
\om_{d-1}(\Th) \ge { d \, \om_d(\Th)  \over \om_d(\Th) +  d + 1}.
\leqno{(2.3)}
$$
Furthermore,
$$
\om_{n-2} (\Th) \ge { ( \omc_{n-1} (\Th)  - 1) \om_{n-1} (\Th)
\over \om_{n-1} (\Th) + \omc_{n-1} (\Th) }. \leqno (2.4)
$$
 
 \medskip
The lower bounds (2.1) and (2.3) are implicit in \cite{Schm67} and are stated 
in \cite{Lau08b}. It is shown
in \cite{Lau08b} that their composition produces Khintchine's theorem. 
The same splitting principle
is used here. We prove Theorem 1 in \S 7    
by iterating successively the finer
Going-up estimates (2.2) and (2.1), and in the 
other direction the Going-down inequalities (2.4) and (2.3). 

In contrast with the previous works \cite{Lau08a,Lau08b,Schm67}, 
our approch is based here on  the use of the second theorem of
Minkowski on the successive minima of a convex body,  
combined with Mahler's theory of compound
convex bodies \cite{Mah}. 

We conclude this section by formulating the transfer inequalities
between $\om_d (\Th)$ and $\om_{d'} (\Th)$ that easily
follow from repeated applications of (2.1) and (2.3).

\proclaim Corollary 1.
Let $\Th  =(\th_1, \dots , \th_n)$ be in $\bR^n$ with 
$1, \th_1, \dots , \th_n$ linearly independent over $\bQ$.
For any integers $d,d'$ with $0\le d < d' \le n-1$, we have
$$
{(d+1) \om_{d'}(\Th)\over (d'-d)\om_{d'}(\Th) + d'+1} 
\le \om_d (\Th) \le {(n-d') \om_{d'}(\Th) - d' + d \over n - d}.
$$

\section{3. Multilinear algebra}
We collect in this section some classical results of 
multilinear algebra and their geometrical interpretation
in terms of join
and intersection of linear varieties in the space $\bR^{n+1}$.
For more details, we refer to \cite{Bou}. 

First, we equip the real vector space $\bR^{n+1}$ with the usual scalar product 
$$
\bx \cdot  \by = x_1y_1+Ê\cdots +  x_{n+1}y_{n+1} , 
\quad \bx =(x_1, \dots , x_{n+1}), \quad \by=(y_1, \dots, y_{n+1}),
$$
and extend it naturally to the Grassmann algebra 
$\Lambda(\bR^{n+1})$, by requiring that for any orthonormal 
basis $\{\be_i\}_{1\le i\le n+1}$ of $\bR^{n+1}$, the family of  wedge products 
$$
\be_{i_1}\wedge \dots \wedge \be_{i_r} \, ;  
\quad  1 \le i_1 < \cdots < i_r \le n+1, \,\, 0 \le r\le n+1, 
$$
is an orthonormal basis of $\Lambda(\bR^{n+1})$.  
Then, the Cauchy-Binet formula shows that
$$
\bX\cdot\bY = \det\Big(\bx_i\cdot\by_j\Big)_{1\le i, j \le r} \leqno{(3.1)}
$$
for any pair of decomposable $r$-vectors 
$\bX= \bx_1\wedge\dots \wedge \bx_r$ and  $\bY= \by_1\wedge\dots \wedge \by_r$.
The scalar product $\cdot$  enables us to identify  
the dual of  the real vector space $\Lambda^r(\bR^{n+1})$ with itself. 
For any multivector
$\bX \in \La(\bR^{n+1})$, we denote by 
$ | \bX | = \sqrt{\bX \cdot \bX}$ the Euclidean norm of $\bX$.

Let $\bX\in \La^r(\bR^{n+1})$ and $\bY\in \La^s(\bR^{n+1})$ 
be two multivectors of respective degree $r$ and $s$ with $s\le r$.
We define the {\it internal product} 
(also called {\it contraction} ) of $\X$ by $\bY$, as the unique multivector
$$
\bY \lrcorner\bX \in \La^{r-s} (\R^{n+1}) 
$$
for which the equality 
$$
\bZ\cdot (\bY\lrcorner\bX) = (\bZ \wedge \bY)\cdot \bX \leqno{(3.2)}
$$
holds for any $ \bZ \in  \La^{r-s} (\R^{n+1})$. 
In other words, the application $\bX\mapsto \bY\lrcorner\bX$
is the transpose of the linear map $\bZ \mapsto \bZ\wedge \bY$ 
with respect to the dot pairing.

Assume  now that  $\bX= \bx_1\wedge\dots \wedge \bx_r$ and  
$\bY= \by_1\wedge\dots \wedge \by_s$  
are decomposable multivectors. 
When $s=1$, we deduce from $(3.1)$ and $(3.2)$ the explicit  formula 
$$
\by\lrcorner\bX = \sum_{j=1}^r (-1)^{r-j} (\by\cdot\bx_j) \bx_1   
\wedge \ldots \wedge {\hat \bx_j} \wedge \ldots \wedge \bx_r  
\leqno{(3.3)} 
$$
for any vector $\by \in \La^1 (\R^{n+1})$. 
It formally follows from (3.2) that 
$$
(\bY\wedge\bY')\lrcorner\bX = \bY\lrcorner(\bY'\lrcorner\bX) \leqno{(3.4)}
$$ 
for any pair of multivectors $\bY$ and $\bY'$ with 
respective degree $s$ and $s'$ such that
$s+s' \le r$. 
Starting with  (3.3) and using (3.4), we obtain by induction on $s$  the formula 
$$
\bY\lrcorner\bX = 
\sum {\rm sgn}(\sigma) (\by_1\cdot \bx_{\sigma(r-s+1)})\cdots 
(\by_s\cdot\bx_{\sigma(r)}) \bx_{\sigma(1)}\wedge\ldots \wedge\bx_{\sigma(r-s)}  
\leqno{(3.5)}
$$
where the sum is taken over all  the substitutions 
$\sigma$ of $\{1,\ldots,r\}$ such that  $\sigma(1)<  \cdots < \sigma(r-s)$.

\medskip

Let $\{\be_1, \ldots , \be_{n+1}\}$ be any  
positively oriented  (meaning that $\det(\be_1, \ldots ,  \be_{n+1})=1$)
 orthonormal  basis of $\R^{n+1}$. Remark that the  
volume form  $\be_1\wedge \ldots \wedge \be_{n+1}$
 does not depend upon the choice of such a basis.

\proclaim Definition 4. For every $\bX$ in $\La^r (\R^{n+1})$, we denote by
$$
\ast {\bX} = \bX \lrcorner(\be_1 \wedge \ldots 
\wedge \be_{n+1})\in \La^{n+1-r}(\bR^{n+1})
$$
the {\it Hodge dual} of $\bX$.

Expanding 
$$
\bX= \sum_{1\le i_1 < \cdots < i_r\le n+1} 
X_{i_1,\ldots ,i_r} \be_{i_1}\wedge \ldots \wedge \be_{i_r}
$$
in the induced orthonormal  basis of $\La^r(\bR^{n+1})$, we find
$$
\ast \bX = \sum_{1\le i_1 < \cdots < i_r\le n+1}
\varepsilon_{i_1,\ldots ,i_r} X_{i_1,\ldots ,i_r} 
\be_{j_1}\wedge \ldots \wedge \be_{j_{n+1-r}}, 
$$
where $\{j_1,\ldots , j_{n+1-r}\} = \{1, \ldots , n+1\} 
\setminus \{i_1, \ldots ,i_{r}\}$
with $j_1< \cdots < j_{n+1-r}$,  and $\varepsilon_{i_1,\ldots ,i_r}$ 
stands for the signature of the shuffle
substitution $(1, \ldots, n+1) \mapsto
( j_1,\ldots , j_{n+1-r}, i_1,\ldots , i_r)$. The Hodge star operator
$$
\ast : \La^{r}(\bR^{n+1}) \buildrel \sim 
\over \longrightarrow  \La^{n+1-r}(\bR^{n+1})
$$
is clearly an isometry for the dot scalar product 
and iterating  twice the Hodge star, we get 
$$
\ast  \circ \ast = (-1)^{r(n+1-r)} {\rm Id}.  \leqno (3.6)  
$$

\proclaim
Lemma 1. 
Let  $\bX = \bx_1 \wedge \ldots \wedge \bx_{r}$ be  a system of Pl\"ucker
coordinates 
\footnote{{\rm $(\natural)$}}{{ \rm The word ``coordinates''  
classically  refers to   the canonical basis of $\La^r(\bR^{n+1})$.}}
of a  $r$-dimensional  subspace 
$$
V = \, < \bx_1, \ldots , \bx_{r} >  
$$
in $\bR^{n+1} $.   Then $\ast{\bX}$ is a system of 
Pl\"ucker coordinates of the  orthogonal $V^{\perp}$ of $V$. 

\pro
That is the assertion of Theorem I of  Chapter VII \S 3  in \cite{HoPe}.
Using the  notion of contraction, we may argue as follows.  
Take any  orthonormal basis
$\{\be_1, \ldots, \be_{r}\}$ of $V$ and extend it to an  orthonormal 
basis $\{\be_1, \ldots, \be_{n+1}\}$ of $\bR^{n+1}$. Then
$$
\bX= \rho  (\be_1 \wedge \ldots \wedge \be_{r})$$
for some non-zero real number $\rho$. Now, it follows from (3.5) that  
$$
\ast \bX = \pm  \rho (\be_{r+1} \wedge \ldots \wedge \be_{n+1}).
$$
\cqfd

\rem
The same argument shows more generally that 
if $\bY = \by_1\wedge \ldots \wedge \by_s$ is a system of
 Pl\"ucker coordinates of an $s$-dimensional vector space
$W = \, < \by_1, \ldots , \by_{s} >$ with $s\ge r$, 
then $\bX\lrcorner \bY$ is a system of Pl\"ucker coordinates of
the intersection $W\cap V^\perp$, provided that 
this intersection   has dimension $s-r$.
\medskip

\proclaim Lemma 2.
For any $\bX \in \La^r (\R^{n+1})$ and $\bY \in \La^s (\R^{n+1})$
with $r + s \le n+1 $, we have the duality formul\ae
$$
\ast({\bY \wedge \bX}) =  \bY \lrcorner({\ast  \bX}) 
$$

\pro Using (3.4), we find 
$$
\ast({\bY \wedge \bX}) = ({\bY \wedge \bX}) \lrcorner
(\be_1\wedge \ldots \wedge \be_{n+1}) 
= \bY\lrcorner(\bX\lrcorner(\be_1\wedge \ldots \wedge \be_{n+1}))
=\bY\lrcorner(\ast \bX).  
$$
\cqfd

\section{4. Alternative definition of the intermediate exponents}
Let $P$ and $Q$ be points in $\bP^n(\bR)$
with homogeneous coordinates ${\bx}$ and ${\by}$. As in \cite{Lau08b}, 
we define the projective distance ${\rm d}(P,Q)$
between $P$ and $Q$  by  
$$
{\rm d}(P,Q) ={ |{\bx} \wedge {\by} | \over | {\bx} | | {\by} |}.
$$
 It has been shown in Lemma 1 of \cite{Lau08b} 
that for any point $\Th$ in
$\bP^n(\bR)$ with homogeneous coordinates $\by=(1, \th_1, \ldots , \th_n)$ and any 
linear subvariety $L$ of  $\bP^n(\bR)$ with Pl\"ucker coordinates $\bX$, the 
minimal distance
${\rm d} ( \Th , L) $ between $\Th$ and the set of 
real points of $L$ is equal to 
$$
{\rm d} ( \Th , L)  = { | \by \wedge \bX | \over | \by | | \bX |}.  \leqno{(4.1)}
$$
We  can now reformulate Definition 3 in terms  of 
integer solutions of the following  system of linear inequations.

\proclaim
Proposition.  For any integer $d$ with $0 \le d \le n-1$, 
the exponent $\om_d(\Th)$ is the supremum 
of the real numbers $\om$
for which there exist infinitely many integer  
multivectors $\bX \in \Lambda^{d+1}(\bZ^{n+1})$ such that
$$
| \by \wedge \bX | \le | \bX |^{-\om}.
$$

In relation with  Definition 4 of  \cite{Lau08b}, we do not 
assume here that the multivectors $\bX$ occurring in the Proposition 
are decomposable. To suppress this additional condition, 
we expand  the remark given on page 312 of \cite{Lau08b}. The following lemma 
will be as well our main ingredient to prove  Theorem 2.

\proclaim
Lemma 3.
Let $\by = (1, \theta_1, \ldots , \theta_n) \in \R^{n+1}$ and let 
$U, V$ be positive real numbers with $V \le U$.
The convex body $\cC$ of $\Lambda^{d+1} (\R^{n+1})$ 
consisting  of the $\bZ$ such that
$$
|\bZ| \le U V^d \and  |\by \wedge \bZ| \le V^{d+1} \leqno{(4.2)}
$$
is comparable \footnote{{\rm(\dag)}}{{\rm 
We say that two families $\cC_1$ and $\cC_2$ of symmetrical convex bodies,  
parametrized by (say)  $U$ and $V$,  are {\it comparable} if  there exists 
 a real number $\kappa > 1$,   
  such that the inclusions $\kappa^{-1} \, \cC_1(U,V) \subseteq \cC_2(U,V)  
\subseteq \kappa \cC_1(U,V)$ hold for any  parameters $U,V$.  
 Accordingly,  the constants implied in the forthcoming symbols $\ll$, 
$\gg$ and $\asymp$ may depend on $n$ and $\Th$, but not on $U$ and $V$. 
The relation $f\asymp g$ means that we have both $f\ll g $ and $f\gg g$.
}} 
 to the $(d+1)$-th compound of the convex body $\cC'$
consisting  of the $\bz \in \R^{n+1}$ such that 
$$
|\bz| \le U \and  |\by \wedge \bz| \le V.  \leqno (4.3) 
$$

\pro
The convex body $\cC'$ 
is comparable to the  
parallelepiped $\cP$ defined by
$$
|x_0| \le U, \quad |x_0 \theta_i - x_i| \le V, \quad 
1 \le i \le n.
$$
However, $\cP$ is the convex hull of the points
$$
\pm U\by, \pm V \be_1, \ldots , \pm V \be_n,
$$
where
$$
\be_1 = (0, 1, 0, \ldots , 0), \ldots , \be_n = (0, \ldots , 0, 1).
$$

The convex compound $\cC'^{d+1}$ is then comparable to
the convex hull in $\Lambda^{d+1} (\R^{n+1})$ of the
exterior products of $d+1$ of these points, that is, of
$$
\pm V^{d+1} \be_{i_0} \wedge \ldots \wedge \be_{i_d}, \quad
1 \le i_0 < \cdots < i_d \le n, 
$$
and
$$
\pm U V^d \by \wedge \be_{i_1} \wedge \ldots \wedge \be_{i_d}, \quad
1 \le i_1 < \cdots < i_d \le n. 
$$
The points $\bZ$ of this form satisfy 
$$
|\bZ| \ll UV^d, \quad |\by \wedge \bZ| \ll  V^{d+1}.
$$

Conversely, let $\bZ$ be in $\Lambda^{d+1} (\R^{n+1})$ 
for which (4.2) holds
and express it
in the base composed of the $d+1$ exterior products of the
base $(\by, \be_1, \ldots , \be_n)$, that is,
$$
\bZ = \sum \, a_{i_0, i_1, \ldots, i_d} \, \be_{i_0} \wedge \ldots \wedge \be_{i_d} 
+ \sum \, b_{i_1, i_2, \ldots, i_d} \,  
\by \wedge \be_{i_1} \wedge \ldots \wedge \be_{i_d}.
$$
Then, we have the estimates  
$$
\sum \, | a_{i_0,  i_1, \ldots ,i_d} |  + \sum \, | b_{i_1, i_2, \ldots, i_d}|   
 \asymp   | \bZ |   \le  UV^d 
\and  
 \sum \, | a_{i_0,i_1, \ldots, i_d} |   \asymp  | \by \wedge \bZ |   \le  V^{d+1}. 
  $$
 This completes the proof of the lemma. \cqfd

\vskip 2mm

With this lemma, we are able to establish our Proposition.

\noindent {\bf Proof of the Proposition. }
Let $\om$ be a real number with $\om \ge -1$ and
let $\bX$ be a non-zero point in $\Lambda^{d+1} (\bZ^{n+1})$  such that
$$
|\by \wedge \bX| \le |\bX|^{-\om}. 
$$
The first minimum  of the convex body $\cC$ composed of the 
$\bZ \in \Lambda^{d+1} (\R^{n+1})$ such that
$$
|\bZ| \le |\bX|  \and 
|\by \wedge \bZ| \le |\bX|^{-\om}
$$
is therefore at most equal to $1$ since $\bX$ belongs to $\cC$.
Setting
$$
U = |\bX|^{(d\om +d+1)/(d+1)}, \quad
V = |\bX|^{-\om /(d+1)},
\leqno{(4.4)}
$$
we observe that $V\le U$ and that 
$$
|\bX| = U V^d, \quad |\bX|^{-\om }=V^{d+1}.
$$
By Lemma 3, the convex $\cC$ is comparable to
the $(d+1)$-th compound of the convex body $\cC' \subset \bR^{n+1}$
defined by the inequalities (4.3).
Now, Mahler's theory on compound convex bodies tells us that the integer point
where $\cC$ reaches its first minimum 
is essentially obtained as the wedge product
$\bx_1\wedge \ldots  \wedge \bx_{d+1}$ of the 
integer points $\bx_i, 1\le i\le d+1$, where 
$\cC'$ reaches its $i$-th minimum. We may therefore 
assume that $\bX= \bx_1\wedge \ldots  \wedge \bx_{d+1}$.
Let $L \subset \bP^n(\bR)$ be the $d$-dimensional rational   
linear subvariety $L=\bP(V)$ where 
$V= \langle \bx_1, \ldots , \bx_{d+1}\rangle$. 
 By (4.1), we obtain
$$
{\rm d} ( \Th , L)  = { | \by \wedge \bX | \over | \by | | \bX |} 
\le | \by|^{-1} | \bX |^{-1-\om}  \ll H(L)^{-1-\om}, 
$$
so that $\om_d(\Th)\ge \om$.

Conversely, if $L$ satisfies ${\rm d} ( \Th , L)  \le H(L)^{-1-\om}$, 
choose a system of coprime integer Pl\"ucker coordinates
$\bX$ of $L$, so that $H(L)=| \bX|$. Then  (4.1) shows 
that the upper bound $|\by \wedge \bX| \ll |\bX|^{-\om}$  holds true.  \cqfd 

\medskip
Our Proposition  enables us to recover the following corollary,  
which was already obtained in \cite{Lau08b} and earlier in \cite{Schm67},
using different arguments.
\proclaim
Corollary 2.
For any integer $d$ with $0 \le d \le n-1$, we have the lower bound 
$$
 \om_d(\Th) \ge {d+1\over n-d}.
 $$

\pro The map $\La^{d+1}(\bR^{n+1}) \longrightarrow  \La^{d+2}(\bR^{n+1})$  
which sends $\bX\mapsto \by\wedge \bX$ has rank
${n+1\choose d+1} -{n\choose d}$. Applying the 
Box Principle to the system of linear forms $\by\wedge\bX$
in the integer variables $\bX$, we find that
$$
\om_d(\Th )\ge  { {n+1 \choose d+1}- \left( {n+1\choose d+1} -{n\choose d} \right)\over {n+1\choose d+1} -{n\choose d}}  
= { {n \choose d}\over {n+1\choose d+1} -{n\choose d}} 
= {d+1\over n-d},
$$
as claimed. \cqfd

\section{5. Proof of Theorem 2}
We use the Proposition  as a more convenient characterization  
of the exponents $\om_d(\Th)$ and 
take again the notations of Section 4.
Let $\om$ be a real number with $-1 \le  \om < \om_d (\Th)$ and
let $\bX \in \Lambda^{d+1} (\Z^{n+1})$ be such that
$$
|\by \wedge \bX| \le |\bX|^{-\om},
$$
where $\by$ denotes the homogeneous coordinates of $\Th$.
 Recall that $U$ and $V$ are given by (4.4) and that
the convex bodies $\cC$ and $\cC'$ are 
defined by (4.2) and (4.3), respectively.
The first minimum $\lambda_1$ of the convex body $\cC$ 
is at most equal to $1$ since $\bX$ belongs to $\cC$.
Replacing possibly $\bX$ by the integer point where 
this first minimum is reached and increasing suitably  $\om$,
 we may   assume that $\lambda_1=1$.

By Lemma 3, the convex $\cC$ is comparable to
the $(d+1)$-th compound of the convex body $\cC'$
of volume
$$
\vol(\cC') \asymp  U V^n = |\bX|^{(-(n-d)\om +d+1)/(d+1)}.
$$
By Minkowski's Theorem, the successive minima
$\lambda'_1 \le \ldots \le \lambda'_{n+1}$ of
$\cC'$ satisfy
$$
\lambda'_1 \times \ldots \times \lambda'_{n+1} \asymp
\vol(\cC')^{-1}  \asymp
|\bX|^{((n-d)\om -d-1)/(d+1)}.
$$
Since $\cC$ is comparable to the $(d+1)$-th
compound of $\cC'$, Mahler's theorem on compound convex bodies
asserts that $\lambda_1$, the first minimum of $\cC$,
is comparable to the product 
$\lambda'_1 \times \ldots \times \lambda'_{d+1}$.
Consequently, 
$$
\lambda'_1 \times \ldots \times \lambda'_{d+1} \asymp 1 \leqno (5.1)
$$
and
$$
(\lambda'_{d+2})^{n-d} \le
\lambda'_{d+2} \times \ldots \times \lambda'_{n+1}
\asymp |\bX|^{((n-d)\om -d-1)/(d+1)},
$$
whence
$$
\lambda'_{d+2} \ll |\bX|^{((n-d)\om-d-1)/((d+1)(n-d))}.  \leqno (5.2)
$$
Now, since the $(d+2)$-th compound of $\cC'$ has its first
minimum comparable to
$$
\lambda'_1 \times \ldots \times \lambda'_{d+2} \asymp \lambda'_{d+2},
$$
it follows from Lemma 3 that
there exists $\bX' \in \Lambda^{d+2} (\Z^{n+1})$ such that
$$
|\bX'| \ll \lambda'_{d+2} U V^{d+1},
\quad  |\by \wedge \bX'| \ll \lambda'_{d+2}  V^{d+2}.
$$
A rapid computation using (5.2) yields that
$$
\lambda'_{d+2} U V^{d+1} \ll |\bX|^{(n-d-1)/(n-d)}
$$
and
$$
\lambda'_{d+2}  V^{d+2} \ll |\bX|^{- ((n-d)\om + 1)/(n-d)}.
$$
This gives
$$
|\by \wedge \bX'| \ll |\bX'|^{- ((n-d)\om + 1)/(n-d-1)},
$$
and we get (2.1) since $\om$ can be taken arbitrarily close to $\om_d (\Th)$.

To establish (2.2), let us first observe that (5.2) with $d=0$ gives
$$
\lambda'_2 \ll |\bX|^{\om - 1/n}. \leqno (5.3)
$$
One can get a better upper bound
for $\lambda'_2$ when $d=0$ by taking the uniform exponents
into account, as we show now. In that case 
$\cC =\cC' $ and $\lambda'_1 =\lambda_1=1$.
 The vector $\bX$ is necessarily primitive in $\bZ^{n+1}$,  since 
the convex body $\cC'$ attains its first minimum at that point.  
Let $\omc$ be a real number 
with $\omc < \omc_0 (\Th)$. 
  By Definition 2, there exists a non-zero integer point $\bx$ such that 
$$
|\bx| < |\bX| , \quad |\by \wedge \bx |  \le |\bX|^{-\omc}.
$$
Since $\bX$ is primitive, the vectors $\bx$ 
and $\bX$ are linearly independent. This gives
$$
\lambda'_2 \ll |\bX|^{\om  -\omc }. \leqno (5.4)
$$
Note that the upper estimate (5.4) may be
sharper than (5.3) since $\omc_0 (\Th) \ge 1/n$. 

Observing  that $U=|\bX|$ and $ V = |\bX|^{-\om}$ and 
proceeding as above, we infer from (5.4) that
$$
\lambda'_2 U V \ll |\bX|^{1 - \omc}
$$
and
$$
\lambda'_2  V^{2} \ll |\bX|^{- (\om + \omc)},
$$
whence
$$
|\by \wedge \bX'| \ll |\bX'|^{- (\om + \omc )/(1 - \omc)}.
$$
Letting $\om$ tends to $\om_0(\Th)$ 
and $\omc$ tends to $\omc_0(\Th)$, this gives
$$
\om_1 (\Th)
\ge {\om_0 (\Th) + \omc_0 (\Th) \over 1 - \omc_0 (\Th)}. 
$$
We have proved (2.2).

\section{6. Proof of Theorem 3}
The proof is parallel to that of Theorem 2. 
We use Hodge duality to reverse the Going-down 
transfer into a Going-up transfer, noting
that the duality permutes  the dimension with the codimension.
 
Let us  start with the  following dual version of the above Proposition. 

\proclaim Lemma 4.
For $d= 0, \ldots , n-1$, the exponent $\om_{d} (\Th)$ 
of a point $\Th$ in $\bR^n$ with homogeneous coordinates 
$\by$ is the
supremum of the real numbers $\om$ 
such that there are infinitely many
$\bX \in \La^{n-d} (\bZ^{n+1})$ with
$$
|\by \lrcorner\bX| \le |\bX|^{-\om}.
$$

\pro
By Lemma 2 and (3.6), we have 
$$
\ast(\by \wedge \ast {\bX}) = (-1)^{(d+1)(n-d)}( \by \lrcorner\bX),
$$
for every $\bX$ in $\La^{n-d} (\R^{n+1})$. 
Note that  $\ast$
maps $\La^{n-d} (\Z^{n+1})$ isometrically onto $\La^{d+1} (\Z^{n+1})$, so that 
$$
 |\ast{\bX}| = |\bX| \and |\by \wedge \ast{\bX}| = |\by \lrcorner\bX|.
$$
Now, replace $\bX$ by $\ast\bX$ in  the Proposition  to conclude. \cqfd

Here is now the dual version of Lemma 3.

\proclaim
Lemma 5.
Let $d$ be an integer with $0\le d\le n-1$ and let 
$U, V$ be positive real numbers with $V \le U$.
The convex body $\cC$ of $\Lambda^{n-d} (\R^{n+1})$ 
consisting  of the $\bZ$ such that 
$$
|\bZ| \le U^{n-d} \and  |\by \lrcorner \bZ| \le U^{n-d-1}V \leqno{(6.1)}
$$
is comparable to the $(n-d)$-th compound of the convex body $\cC'$
composed of the $\bz \in \R^{n+1}$ such that
$$
|\bz| \le U \and |\by \cdot \bz| \le V.  \leqno (6.2)
$$

\pro
Let $\{\be_1, \ldots , \be_n\}$ be an orthonormal 
basis of the orthogonal of $\by$ in $\bR^{n+1}$.
The convex body $\cC'$   is comparable to  the
parallelepiped $\cP$  consisting of the points 
$$
x_0 \by + x_1\be_1+\cdots + x_n\be_n \quad {\rm where}\quad 
|x_0| \le V, \quad |x_i| \le U, \quad 
1 \le i \le n.
$$
Note that  $\cP$ is comparable to the  convex hull of the points 
$$
\pm V\by, \pm U \be_1, \ldots , \pm U \be_n.
$$
The  compound convex body $\cC'^{n-d}$ is then comparable to
the convex hull in $\Lambda^{n-d} (\R^{n+1})$ of the
exterior products of $n-d$ of these points, that is, of
$$
\pm U^{n-d} \be_{i_1} \wedge \ldots \wedge \be_{i_{n-d}}, \quad
1 \le i_1 < \cdots < i_{n-d} \le n,\leqno{(6.3)}
$$
and
$$
\pm U^{n-d-1} V  \be_{i_1} \wedge \ldots \wedge \be_{i_{n-d-1}}\wedge\by, \quad
1 \le i_1 < \cdots < i_{n-d-1} \le n. \leqno{(6.4)}
$$
Express now any point  $\bZ$  in $\Lambda^{n-d} (\R^{n+1})$ 
in the base composed of the $n-d$ exterior products of the
base $( \be_1, \ldots , \be_n, \by)$, that is,
$$
\bZ = \sum \, a_{i_1,   \ldots,  i_{n-d}} \, \be_{i_1} \wedge \dots \wedge \be_{i_{n-d}}
+ \sum \, b_{i_1,   \ldots , i_{n-d-1}} \, 
  \be_{i_1} \wedge \dots \wedge \be_{i_{n-d-1}}\wedge\by.
$$
Then, formula (3.3) shows that
$$
\by\lrcorner\bZ = |\by |^2 \left(\sum \, b_{i_1,  \ldots,  i_{n-d-1}} \, 
  \be_{i_1} \wedge \dots \wedge \be_{i_{n-d-1}}\right) .
  $$
Therefore, the points  $\bZ$ of the form (6.3) or (6.4) satisfy
$$
|\bZ| \le U^{n-d} \and  |\by \lrcorner \bZ| \le |\by|^2  U^{n-d-1}V. \leqno{(6.5)} 
$$
Conversely, for any point $\bZ$ satisfying (6.5),  
the coefficients $a_{i_1, \ldots,  i_{n-d}}$ (resp. $b_{i_1,  \ldots,  i_{n-d-1}}$) are
bounded in absolute value by $U^{n-d}$ (resp. by $U^{n-d-1}V$).
This completes the proof of the lemma. \cqfd

\vskip 2mm

With these two lemmata, we are able to establish Theorem 3.

\noindent {\bf Proof of Theorem 3.}
Let $\om$ be a  real number with $-1\le \om < \om_{d} (\Th)$. By Lemma 4, 
 there exist infinitely many  points $\bX \in \Lambda^{n-d} (\Z^{n+1})$ such that
$$
|\by \lrcorner\bX| \le |\bX|^{-\om}.
$$
Fix such a point $\bX$ with large norm $|\bX|$ 
and consider the convex body $\cC$ composed of the multivectors 
$\bZ \in \Lambda^{n-d} (\R^{n+1})$ such that
$$
|\bZ| \le |\bX| \and  |\by \lrcorner\bZ| \le |\bX|^{-\om}.
$$
It contains the integer point $\bX$. Replacing possibly $\bX$ by  
a smaller point and enlarging suitably $\om$,  one can
assume that $\bX$ is the smallest non-zero integer point in $\cC$. 
Thus, we may assume that the  first minimum of $\cC$  is equal to $1$.
Setting
$$
U = |\bX|^{1/(n-d)} \and 
V = |\bX|^{-((n-d) \om +n-d - 1)/(n-d)},
$$
we observe that $V\le U$ and that 
$$
|\bX| = U^{n-d}, \quad |\bX|^{-\om}=U^{n-d-1} V.
$$
By Lemma 5,  the convex body $\cC$ is therefore 
comparable to the $(n-d)$-th compound of the convex body $\cC'$ 
consisting  of the real $(n+1)$-tuples $\bz$  such that  
$$
|\bz| \le U \and  |\by \cdot \bz| \le V.
$$
Let 
$$
\lambda_1 \le \ldots \le \lambda_{n+1}
$$
be the successive minima of the convex body $\cC'$. 
Since the  Euclidean volume of $\cC'$  is $\asymp U^nV$,
the second theorem of Minkowski gives
$$
\lambda_1 \times \ldots \times \lambda_{n+1} \asymp 
(U^n V)^{-1}= | \bX |^{((n-d)\om -d-1)/(n-d)}.
$$
Since the first minimum of the $(n-d)$-th compound of $\cC'$
is comparable to $1$, one gets
$$
\lambda_1 \times \ldots \times \lambda_{n-d} \asymp 1,
$$
hence
$$
\lambda_{n-d+1 } \times \ldots \times \lambda_{n+1} 
\asymp | \bX |^{((n-d)\om -d-1)/(n-d)}.
$$
Consequently,
$$
\lambda_{n-d+1}^{d+1} \ll  | \bX |^{((n-d)\om -d-1)/(n-d)},\leqno{(6.6)}
$$
and
$$
\lambda_{n-d+1} U \ll |\bX|^{\om/(d+1)}.
$$
Since the first minimum of the $(n-d+1)$-th
compound of $\cC'$ is comparable to
the product $\lambda_1 \times \ldots \times \lambda_{n-d+1}$,
hence to $\lambda_{n-d+1}$, we infer from Lemma 5
that there exists a non-zero
integer point $\bX' \in \Lambda^{n-d+1} (\Z^{n+1})$
such that
$$
|\bX'| \ll \lambda_{n-d+1} U^{n-d+1} = \lambda_{n-d+1}U|\bX| \ll |\bX|^{(\om+d+1)/(d+1)}
$$
and
$$
|\by \lrcorner\bX'| \ll \lambda_{n-d+1} U^{n-d} V = \lambda_{n-d+1} U | \bX|^{-\om} 
\ll |\bX|^{-d\om/(d+1)}.
$$
Since $\om$ can be taken arbitrarily close to $\om_d(\Th)$, Lemma 4 gives  (2.3).

For $d=n-1$, it is possible to get a sharper result.
 In that case $\cC =\cC' $ is a convex body  in $\bR^{n+1}$ and  (6.6) reads
 $$
 \lambda_2 \ll |\bX|^{-1+\om/n}. 
 \leqno{(6.7)}
 $$
 Enlarging possibly $\om$, we may assume that 
 $$
  | \by \cdot \bX | = | \bX|^{-\om} .
 $$
 The vector $\bX$ is necessarily primitive in $\bZ^{n+1}$,  since 
the convex body $\cC'$ attains its first minimum at that point.  
Let $\omc$ be a real number 
with $\omc < \omc_{n-1} (\Th)$. 
  By Definition 2, there exists a non-zero integer point $\bx \in \bZ^{n+1}$ such that 
$$
|\bx| \le  |\bX|^{\om/\omc}  \and  |\by \cdot \bx |  < |\bX|^{-\om}.
$$
Since $\bX$ is primitive, the vectors $\bx$ and $\bX$ 
are linearly independent; otherwise $\bx$ should be an integer multiple
of $\bX$ and $ |\by \cdot \bx |$ should be greater than or equal to 
$ | \by \cdot \bX | = | \bX|^{-\om}$. Thus, we obtain the upper bound
$$
\lambda_2 \ll |\bX|^{-1+\om/\omc },  \leqno{(6.8)}
$$
which may  be better than (6.7) since $\omc_{n-1}(\Th) \ge  n$. 
Now, we take again the preceding arguments. 
Noting that $U= |  \bX | $ and $ V=|\bX|^{-\om}$,
 we obtain a non-zero point $\bX'\in \La^2(\bZ^{n+1})$ satisfying
$$
|\bX'| \ll \lambda_{2} U^{2}  \ll |\bX|^{1 +  \om/\omc}
$$
and
$$
|\by \lrcorner\bX'| \ll \lambda_2 U V \ll |\bX|^{-\om+\om/\omc}.
$$
Then, Lemma 4  gives
$$
\om_{n-2} (\Th) \ge {(\omc - 1) \om 
\over \om + \omc }.
$$
Letting $\om$ and $\omc$ tend respectively to $\om_{n-1}(\Th)$ 
and $\omc_{n-1}(\Th)$, we have established (2.4).

\section{7. Proof of Theorem 1}
It is a formal consequence of the finer estimates (2.1)--(2.4).

Using the second  inequality of Corollary 1 with  $d=1$ and $d'=n-1$, we get the estimate 
$$
\om_{n-1} (\Th) \ge (n-1) \om_1 (\Th) + n - 2,
$$
which, combined with (2.2), yields the right 
hand side of the claimed  estimate, namely
$$
\om_{n-1} (\Th) \ge (n-1) 
{\om_0 (\Th) + \omc_0 (\Th) \over 1 - \omc_0 (\Th)} + n - 2
={(n-1)\om_0(\Th) +\omc_0 (\Th)+n-2\over 1 - \omc_0 (\Th)}.
$$

Using now the first inequality of Corollary 1 with $d=0$ and $d' =n-2$,  we get  
$$
\om_0 (\Th) \ge {\om_{n-2} (\Th) \over (n - 2) \om_{n-2} (\Th)
+n - 1}
$$
which, combined with (2.4), yields the claimed Going-down transfer
inequality, namely
$$
\om_0 (\Th) \ge {(\omc_{n-1}(\Th) -1)\om_{n-1}(\Th) \over 
((n-2)\omc_{n-1}(\Th) +1)\om_{n-1}(\Th) +(n-1)\omc_{n-1}(\Th)}.
$$

\vskip 8mm

\bigskip

\centerline{\bf References}

\vskip 5mm

\beginthebibliography{999}

\bibitem{Bou}
N. Bourbaki,  Algebra 1, Chapter 3, Springer-Verlag, New-York, 1989. 

\bibitem{BuLa05a}
Y. Bugeaud and M. Laurent,
{\it Exponents of Diophantine approximation and Sturmian
continued fractions},
Ann. Inst. Fourier (Grenoble) 55 (2005), 773--804.

\bibitem{BuLa05b}
Y. Bugeaud and M. Laurent,
{\it Exponents of inhomogeneous Diophantine approximation},
Moscow Math. J. 5 (2005), 747--766.

\bibitem{BuLa07}
Y. Bugeaud and M. Laurent,
{\it Exponents of Diophantine approximation}.
In: Diophantine Geometry Proceedings,
Scuola Normale Superiore Pisa, Ser. CRM, vol. 4, 2007, 101--121.

\bibitem{BuLiv}
Y. Bugeaud,
Approximation by algebraic numbers.
Cambridge Tracts in Mathematics, Cambridge, 2004.

\bibitem{Cas}
J. W. S. Cassels,
An introduction to Diophantine Approximation.
Cambridge Tracts in Math. and Math. Phys., vol. 99, Cambridge
University Press, 1957.

\bibitem{Dys}
F. J. Dyson,
{\it On simultaneous Diophantine approximations},
Proc. London Math. Soc. 49 (1947), 409--420.

\bibitem{GrLe}
P. M. Gruber and C. G. Lekkerkerker,
Geometry of numbers.
Series Bibliotheca Mathematica 8,
North--Holland, Amsterdam, 1987.

\bibitem{HoPe}
H. Hodge and D. Pedoe, Methods of algebraic  Geometry, Cambridge University Press, 1947.

\bibitem{Jar35}
V. Jarn\'\i k, 
{\it \"Uber einen Satz von A. Khintchine},
Pr\'ace Mat.-Fiz. 43 (1935), 1--16.

\bibitem{Jar36}
V. Jarn\'\i k, 
{\it \"Uber einen Satz von A. Khintchine, 2. Mitteilung},
Acta Arith. 2 (1936), 1--22.

\bibitem{Kh26}
A. Ya. Khintchine,
{\it \"Uber eine Klasse linearer diophantischer Approximationen},
Rendiconti Circ. Mat. Palermo 50 (1926), 170--195.

\bibitem{Lau08a}
M. Laurent,
{\it Exponents of Diophantine Approximation in dimension two},
Canad. J. Math. To appear in 2008.

\bibitem{Lau08b}
M. Laurent,  
{\it On transfer inequalities in Diophantine Approximation}.
In: ``Analytic Number Theory in Honour of Klaus Roth''. 
Cambridge University Press (to appear in 2009), 306--314. 

\bibitem{Mah36}
K. Mahler,
{\it Neuer Beweis einer Satz von A. Khintchine},
Mat. Sbornik 43 (1936), 961--962.

\bibitem{Mah}
K. Mahler,
{\it On compound convex bodies, I},
Proc. London Math. Soc. (3)  {5} (1955), 358--379.

\bibitem{Schm67}
W. M. Schmidt,
{ \it On heights of algebraic subspaces and diophantine
approximations},
Annals of Math. 85 (1967), 430--472.

\bibitem{SchmLN}
{W. M. Schmidt},
Diophantine Approximation. Lecture Notes in Math. 
{785}, Springer, Berlin, 1980.

\endthebibliography

\vskip1cm

\noindent Yann Bugeaud  \hfill{Michel Laurent}

\noindent Universit\'e Louis Pasteur
\hfill{Institut de Math\'ematiques de Luminy}

\noindent U. F. R. de math\'ematiques
\hfill{C.N.R.S. -  U.M.R. 6206 - case 907}

\noindent 7, rue Ren\'e Descartes      \hfill{163, avenue de Luminy}

\noindent 67084 STRASBOURG  (FRANCE)
\hfill{13288 MARSEILLE CEDEX 9  (FRANCE)}

\vskip2mm

\noindent {\tt bugeaud@math.u-strasbg.fr}
\hfill{{\tt laurent@iml.univ-mrs.fr}}

\bye